\input amstex\documentstyle{amsppt}  
\pagewidth{12.5cm}\pageheight{19cm}\magnification\magstep1  
\topmatter
\title Open problems on Iwahori-Hecke algebras\endtitle
\author G. Lusztig\endauthor
\address{Department of Mathematics, M.I.T., Cambridge, MA 02139}\endaddress
\thanks{Supported by NSF grant DMS-1855773.}\endthanks
\endtopmatter   
\document

\define\lb{\linebreak}

\define\part{\partial}

\define\n{\notin}
\define\iy{\infty}
\define\m{\mapsto}

\define\sub{\subset}    

\define\T{\times}

\define\nl{\newline}
\redefine\i{^{-1}}

\define\un{\underline}
\define\ov{\overline}
\define\ot{\otimes}

\define\g{\gamma}

\define\Si{\Sigma}

\define\CC{\bold C}

\define\NN{\bold N}

\define\ZZ{\bold Z}

Below we state four open problems (see (a)-(d)) on Iwahori-Hecke algebras.

\subhead 1\endsubhead
Let $I$ be a finite set and let $(m_{ij})_{(i,j)\in I\T I}$ be a symmetric matrix whose diagonal entries are $1$
and whose nondiagonal entries are integers $\ge2$ or $\iy$. Let $W$ be the group with generators
$\{s_i;i\in I\}$ and relations $(s_is_j)^{m_{ij}}=1$ for any $i,j$ such that $m_{ij}<\iy$; this is a Coxeter
group. (Examples of Coxeter groups are the Weyl groups of simple Lie algebras; these are finite groups. Other
examples are the affine Weyl groups which are almost finite.)
For $w\in W$ let $|w|$ be the smallest integer $n\ge0$ such that $w$ is a product of $n$ generators
$s_i,i\in I$. We assume that we are given a weight function $L:W@>>>\NN$ that is, a function such that
$L(w)>0$ for all $w\in W-\{1\}$ and $L(ww')=L(w)+L(w')$ for any $w,w'$ in $W$ such that $|ww'|=|w|+|w'|$.
(For example, $w\m|w|$ is a weight function.)
Let $A=\ZZ[v,v\i]$ where $v$ is an indeterminate. Let $H$ be the free $A$-module with basis $\{T_w;w\in W\}$.
There is a unique structure of associative $A$-algebra on $H$ for which
$(T_{s_i}+v^{-L(s_i)})(T_{s_i}-v^{L(s_i)})=0$ for $i\in I$ and $T_wT_{w'}=T_{ww'}$ for any $w,w'$ in $W$ such
that $|ww'|=|w|+|w'|$; this is the Iwahori-Hecke algebra associated to $W,L$.

For $c\in\CC-\{0\}$ let $H_c=\CC\ot AH$ where $\CC$ is viewed as an $A$-algebra via the ring homomorphism
$A@>>>\CC,v\m c$. Now $H_c$ is also referred to as an Iwahori-Hecke algebra.

(a) {\it Show that the algebras associated in \cite{S20} to a supercuspidal representation of a parabolic
subgroup of a $p$-adic reductive group are (up to extension by a group algebra of a small finite group) of
the form $H_q$ where $q$ is a power of $p$, with $H$ associated to an affine Weyl group $W$ and with $L$
in the collection $\Si_W$ of weight functions on $W$ described in \cite{L91, \S17}, \cite{L95}, \cite{L02}.}
\nl
For example, if $W$ is of affine type $F_4$, then $\Si_W$ consists of all $L$ whose values on $\{s_i;i\in I\}$
are $(1,1,1,1,1)$ or $(1,1,1,2,2)$ or $(2,2,2,1,1)$ or $(1,1,1,4,4)$; if $W$ is of affine type $G_2$, then
$\Si_W$ consists of all $L$ whose values on $\{s_i;i\in I\}$ are $(1,1,1)$ or $(1,1,3)$ or $(3,3,1)$ or
$(1,1,9)$.

The statement analogous to (a) for groups with connected centre over a finite field $F_q$ instead of $p$-adic
groups is known to hold, without the words in parenthesis; in that case, $W$ is a Weyl group and $\Si_W$
consists of the weight functions on $W$ described in \cite{L78, p.35}.

\subhead 2\endsubhead
There is a unique group homomorphism $\bar{}:H@>>>H$ such that $\ov{v^nT_w}T_{w\i}=v^{-n}$ for
$n\in\ZZ,w\in W$; it is a ring isomorphism. Let $H_{\le0}=\sum_{w\in W}\ZZ[v\i]T_w\sub H$.
For any $w\in W$ there is a unique element $c_w\in H_{\le0}$ such that $\bar c_w=c_w$ and
$c_w-T_w\in v\i H_{\le0}$ (see \cite{KL},\cite{L03}). Then $\{c_w;w\in W\}$ is an $A$-basis of $H$.
For $x,y,z$ in $W$ we define $f_{x,y,z}\in A$, $h_{x,y,z}\in A$ by $T_xT_y=\sum_{z\in W}f_{x,y,z}T_z$,
$c_xc_y=\sum_{z\in W}h_{x,y,z}c_z$.

(b) {\it Show that there exists an integer $N\ge0$ such that for any $x,y,z$ in $W$ we have
$v^{-N}f_{x,y,z}\in\ZZ[v\i]$.}
\nl
(See \cite{L03, 13.4}.) If $W$ is finite this is obvious. If $W$ is an affine Weyl group, this is known.

We will now assume that (b) holds. 
With $N$ as in (b), we see that $v^{-N}h_{x,y,z}\in\ZZ[v\i]$ for any $x,y,z$ in $W$. It follows that
for any $z\in W$ there is a unique integer $a(z)\ge0$ such that $h_{x,y,z}\in v^{a(z)}\ZZ[v\i]$
for all $x,y$ in $W$ and $h_{x,y,z}\n v^{a(z)-1}\ZZ[v\i]$ for some $x,y$ in $W$. Hence for $x,y,z$ in $W$
there is a well defined integer $\g_{x,y,z\i}$
such that $h_{x,y,z}=\g_{x,y,z\i}v^{a(z)}\mod v^{a(z)-1}\ZZ[v\i]$.
Let $J$ be the free abelian group with basis $\{t_w;w\in W\}$. For $x,y$ in $W$ we set
$t_xt_y=\sum_{z\in W}\g_{x,y,z\i}t_z$. (This is a finite sum.)

(c) {\it Show that this defines an (associative) ring structure on $J$ (without $1$ in general).}
\nl
Assume now that $W$ is a Weyl group or an affine Weyl group and $L=||$. In this case, (c) is known to be true
and the ring $J$ does have a unit element.

More generally, assume that $W$ is an affine Weyl group and
$L\in\Si_W$ (see (a)); in this case there is a (conjectural) geometric description \cite{L16, 3.11} of the
elements $c_w$. From this one should be able to deduce (c) as well as the well-definedness of the
$\CC$-algebra homomorphism $H_q@>>>\un J$ in \cite{L03, 18.9}, where $H_q$ is as in (a) and $\un J=\CC\ot J$
is independent of $q$.
One should expect that the irreducible (finite dimensional) $\un J$-modules, when viewed as $H_q$-modules,
form a basis of the Grothendieck group of $H_q$-modules. (This is indeed so if $L=||$.) This should provide a
construction of the ``standard modules'' of $H_q$ which, unlike the construction in \cite{L95},\cite{L02},
does not involve the geometry of the dual group.

\subhead 3\endsubhead
Assume that $W$ is finite and that $L=||$. Let $C$ be a conjugacy class in $W$; let
$C_{min}$ be the set of all $w\in C$ such that $|w|$ is minimal. For $w\in C$ let $N^w\in A$
be the trace of the
$A$-linear map $H@>>>H$, $h\m v^{2|w|}T_whT_{w\i}$. We have $N^w\in\ZZ[v^2]$. (Note that $N^w|_{v=1}$ is the
order of the centralizer of $w$ in $W$.) From \cite{GP} one can deduce that for $w\in C_{min}$, $N^w$ depends
only on $C$, not on $w$. We say that $C$ is {\it positive} if $N^w\in\NN[v^2]$.
For example, if $C$ is an elliptic regular conjugacy classes (in the sense of \cite{S74}) then
$C$ is positive (see \cite{L18}). If $W$ is of type $A_n$, the positive conjugacy classes are $\{1\}$ and the
class of the Coxeter element. In the case where $W$ is a Weyl group of exceptional type a complete list of
positive conjugacy classes in $W$ is given in \cite{L18}.

(d) {\it Make a list of all positive conjugacy classes in $W$ assuming that $W$ is a Weyl group of type $B_n$
or $D_n$.}

\enddocument